\documentclass[11pt,a4paper]{article}
\usepackage[margin=26mm]{geometry}
\usepackage[T1]{fontenc}
\usepackage{lmodern,microtype}
\usepackage{amsmath,amssymb,amsthm,mathtools,booktabs}
\usepackage[colorlinks=true,linkcolor=blue,citecolor=blue,urlcolor=blue]{hyperref}
\usepackage{fancyhdr}
\allowdisplaybreaks[2]
\newtheorem{theorem}{Theorem}
\newtheorem{lemma}[theorem]{Lemma}

\newtheorem{corollary}[theorem]{Corollary}
\theoremstyle{remark}
\newcommand{\E}{\mathbb E}
\newcommand{\PP}{\mathbb P}
\newcommand{\R}{\mathbb R}
\newcommand{\TV}{\operatorname{TV}}
\newcommand{\Cov}{\operatorname{Cov}}

\title{Proof of the conjecture about inequality of $1/H$-variations of fractional Brownian motion along the equidistant and the Lebesgue partitions}
\author{W. M. Bednorz and R. M. \L ochowski}
\date{Revised 28 October 2026}
\begin{document}
\maketitle
\begin{abstract}
We give an analytic proof of the comparison conjectured by Das,
\L{}ochowski, Matsuda and Perkowski for the fractional Brownian
level-crossing constant: for standard fractional Brownian motion,
$c_H>\E|Z|^{1/H}$ when $0<H<1/2$ and $c_H<\E|Z|^{1/H}$ when
$1/2<H<1$, where $Z\sim N(0,1)$. At $H=1/2$, both constants equal one.
The main tool is Gaussian comparison under ordering of
increment variances. It applies because truncated variation is
translation invariant and lattice-submodular, as follows from its
minimum-total-variation characterization. A tangent variogram produces
Brownian motion with an independent Gaussian drift on the upper side,
and a periodic Brownian bridge on the lower side. The former has an
explicit stationary reflection rate. For the latter we calculate an
integrated cyclic crossing count by first-passage Laplace transforms.
Elementary hyperbolic, theta-series, and Gamma inequalities then give
strict comparisons on both full parameter ranges. No numerical
integration, simulation, interval certificate, or quantitative local
expansion at $H=1/2$ is used.
\end{abstract}

\section{Introduction and statement}
The level-crossing constant of fractional Brownian motion measures the
asymptotic frequency of successive crossings of a fine spatial grid.
The question considered here is how this constant compares with the
Gaussian moment governing power variation along equally spaced time
partitions. Let $B^H$ be standard fractional Brownian motion (fBm),
normalized by
\[
 B^H_0=0,\qquad \E(B^H_t-B^H_s)^2=|t-s|^{2H},\qquad 0<H<1.
\]
For a real continuous path $f$ and $d>0$, write
\[
 \TV^d(f,[0,L])=
 \sup_{0=t_0<\cdots<t_n=L}
 \sum_{i=1}^n (|f(t_i)-f(t_{i-1})|-d)_+.
\]
Allowing partitions whose endpoints lie inside $[0,L]$ gives the same
supremum. The established long-time characterization of the crossing
constant is
\begin{equation}\label{eq:crossing-limit}
 c_H=\lim_{L\to\infty}\frac1L\E\TV^1(B^H,[0,L]);
\end{equation}
Indeed, \cite{BL2025}, Corollary 2.1, gives the finite-time brackets
\[
 \frac{\E\TV^1(B^H,[0,n])}{n}
 \le c_H\le
 \frac{\E\TV^1(B^H,[0,n])+1}{n},\qquad n\in\mathbb N.
\]
Monotonicity in the time horizon extends the resulting limit from integer
times to arbitrary $L\to\infty$. Corollary 3.1 of the same paper
identifies this $c_H$ with the level-crossing constant of \cite{DLMP},
whose strict-sign conjecture is the statement below.
Set
\[
 m_H=\E|Z|^{1/H}
 =\frac{2^{1/(2H)}\Gamma((1+1/H)/2)}{\sqrt\pi},
 \qquad Z\sim N(0,1),\qquad \mu=\E|Z|=\sqrt{2/\pi}.
\]
We write $\psi=(\log\Gamma)'$ and $\psi_1=\psi'$ for the digamma and
trigamma functions, and $\gamma$ for Euler's constant. The series and
integral representations used below are recorded in \cite{DLMF}.

\begin{theorem}\label{thm:main}
For every standard fractional Brownian motion,
\[
 \boxed{
 c_H>m_H\quad(0<H<1/2),\qquad
 c_{1/2}=m_{1/2}=1,\qquad
 c_H<m_H\quad(1/2<H<1).}
\]
\end{theorem}

Thus the Brownian index is the unique parameter in $(0,1)$ at which
the two constants agree. The result determines the sign of their
difference; it does not give an exact formula for $c_H$ away from
$H=1/2$.

The inequalities are unchanged under a different variance normalization:
if $\widetilde B^H=\sigma B^H$ with $\sigma>0$, both $c_H$ and
$\E|B^H_1|^{1/H}$ are
multiplied by $\sigma^{1/H}$. This follows from self-similarity and
\eqref{eq:crossing-limit}, since
$\sigma B^H_t$ has the same law as $B^H_{\sigma^{1/H}t}$.

\subsection{Historical context and earlier results}
The immediate antecedent is the work of Das, \L{}ochowski, Matsuda and
Perkowski \cite{DLMP}, first posted in 2023 and published in 2025.
They established the limit of normalized level-crossing counts and the
corresponding $1/H$-variation along Lebesgue partitions, which are
defined by successive hits of a spatial grid. Their numerical evidence
led to Conjecture~1.6 in Section~1.2 of the arXiv version: the constant
is larger than the equally spaced power-variation constant for $H<1/2$
and smaller for $H>1/2$. Theorem~\ref{thm:main} proves these comparisons
for all parameters on the respective sides of $1/2$.

An earlier conjecture of Cont and Perkowski \cite{CP2019}, appearing in
their 2018 preprint and published in 2019, predicted the normalization
$m_H/2$ for the local-time limit of dyadic upcrossing counts. The
conjecture after Lemma~3.5 in the second arXiv version is a stronger
statement about almost-sure uniform convergence in time and space;
the subsequent paragraph explicitly predicts the coefficient $m_H$ for
$1/H$-variation along dyadic Lebesgue partitions. Thus equality of the
two variation constants was part of that earlier prediction. The
strict inequalities above show that this normalization differs from
the level-crossing normalization whenever $H\ne1/2$.

The truncated-variation formulation seems to make the comparison accessible.
Bednorz and \L{}ochowski \cite{BL2015} developed integrability and
concentration estimates for truncated variation, including fBm.
The preprint \cite{BL2025}, posted in 2025, gives the
finite-time brackets quoted above and identifies the resulting
constant with the crossing constant of \cite{DLMP}. These results
supply the limiting object used here; the sign comparison is proved
in the present note.

The deterministic input is variational characterization
of truncated variation \cite{LochowskiTube}: the truncated variation
at level $d$ is the least total variation of a path within uniform
distance $d/2$. Equation~(1.2) of that paper states the identity, whose
proof combines Corollary~3.8 with Theorem~4.1; it is also recorded as
equation~(1.2) in \cite{BL2015}. The feature needed below is that the
approximating path has no prescribed endpoint values. This freedom
allows the lattice operations used in the Gaussian comparison.

Explicit formulas for Brownian truncated variation were developed in  \cite{LochowskiBM}. In particular, the long-time rate for
Brownian motion with drift is due to \L{}ochowski and Mi\l{}o\'s
\cite{LochowskiMilos}, Fact~11 and equation~(1.7) of the arXiv version.
After variance scaling, this rate is $b\coth(b/v)$ for
$\sqrt v\,W_t+bt$, with value $v$ at $b=0$.
Section~\ref{sec:upper} gives a direct stationary-reflection argument
for the finite-horizon upper bound that the proof requires.

\subsection{The comparison mechanism and proof structure}
Truncated variation optimizes over partitions, so the required estimate
must retain information about the joint Gaussian path. The main
mechanism is that the minimum-variation identity makes the finite-path
functional translation invariant and lattice-submodular. These two
properties turn an ordering of all increment variances into an ordering
of expected truncated variations. Section~\ref{sec:comparison} proves
this statement first for finite vectors and then for continuous paths.

The proof is organized around two explicit bounds. For $H>1/2$, setting
$r=2H-1$, $v=(1-r)T^r$, and $w=rT^{r-1}$ gives
\[
 c_H\le\E\left[\sqrt w\,Z\coth(\sqrt w\,Z/v)\right],\qquad T>0.
\]
For $H<1/2$, setting $p=2H$, $v=(2-p)T^{p-1}$, and
$L=(2-p)T/(1-p)$ gives
\[
 c_H\ge\frac4L\sum_{n\ge1}n e^{-2n^2/(vL)},\qquad T>0.
\]
The upper comparison uses a Brownian motion with an independent
Gaussian drift. On the lower side, the same Brownian bridge is repeated
periodically, and its integrated cyclic crossing count gives a lower
bound valid over arbitrarily many periods. In both cases a tangent
inequality for the fBm variogram selects the comparator.

The remaining scalar estimates use two regimes on each side of the
Brownian parameter. On $1/2<H\le3/4$, a quadratic bound on
$x\coth x$ reduces the upper estimate to a Gamma-function comparison;
on $3/4\le H<1$, a direct Gaussian integral bound gives $c_H<\E|Z|$.
On $1/3\le H<1/2$, a uniform theta-sum estimate controls the full bridge
series; on $0<H\le1/3$, its first term suffices after optimizing the
scale. The estimates overlap at their stated endpoints and require no
separately assumed local theorem near $H=1/2$. All probabilistic
comparisons and scalar estimates used to establish the strict signs
are proved below.

\section{Gaussian increment comparison}\label{sec:comparison}
The free-endpoint tube identity for a c\`adl\`ag path (in particular,
for a continuous path) is
\begin{equation}\label{eq:tube}
 \TV^d(f,[0,A])
 =\inf_{\|u-f\|_\infty\le d/2}\TV(u,[0,A]);
\end{equation}
see \cite{LochowskiTube}, equation (1.2), and \cite{BL2015},
equation (1.2). The infimum ranges over functions
whose endpoint values are subject only to the same tube constraint.
For a finite vector $x=(x_0,\ldots,x_n)$, define truncated variation by
maximizing over increasing subsequences:
\[
 F_d(x)=\max_{0=i_0<\cdots<i_k=n}
       \sum_{j=1}^k (|x_{i_j}-x_{i_{j-1}}|-d)_+.
\]
This is generally different from the sum over adjacent indices. The
minimum-variation identity of \cite{BL2015}, equation (1.2), gives
\[
 F_d(x)=\min\left\{\sum_{i=1}^n|u_i-u_{i-1}|:
                  |u_i-x_i|\le d/2\text{ for every }i\right\}.
\]
To obtain this finite version from the cited path identity, embed $x$ as
the step function $f(t)=x_{\lfloor t\rfloor}$ on $[0,n]$. Sampling a
path approximant at the integers cannot increase its variation, and a
feasible finite vector extends to a step-function approximant of exactly
the displayed variation. The finite minimum is attained by compactness.
The initial and terminal values of $u$ are free subject to the tube
constraint. In particular, this is not an anchored approximation cost.

\begin{lemma}\label{lem:gaussian}
The functional $F_d$ is translation invariant and submodular:
\[
 F_d(x+t\mathbf1)=F_d(x),\qquad
 F_d(x\wedge y)+F_d(x\vee y)\le F_d(x)+F_d(y).
\]
Consequently, if $X,Y$ are centered Gaussian vectors, possibly degenerate, and
\[
 \E(X_i-X_j)^2\le\E(Y_i-Y_j)^2\qquad\text{for all }i,j,
\]
then $\E F_d(X)\le\E F_d(Y)$.
\end{lemma}
\begin{proof}
Translation invariance is immediate. The edge inequality
\[
 |\min(a,c)-\min(b,d)|+|\max(a,c)-\max(b,d)|
 \le |a-b|+|c-d|
\]
implies the lattice inequality for discrete total variation. If $u,v$
are feasible approximants to $x,y$, respectively, then $u\wedge v$ and
$u\vee v$ are feasible approximants to $x\wedge y,x\vee y$. Take
infima to obtain submodularity.

Every coordinate occurs at most twice in any subsequence sum, so
\[
 |F_d(x)-F_d(y)|\le 2\|x-y\|_1.
\]
Let $G$ be a standard Gaussian vector and set
$F_{d,\varepsilon}(x)=\E F_d(x+\varepsilon G)$. For fixed
$\varepsilon>0$, this function is smooth, has bounded first and second
derivatives and at most linear growth, and satisfies
\[
 \sup_x|F_{d,\varepsilon}(x)-F_d(x)|
 \le2\varepsilon\E\|G\|_1\longrightarrow0.
\]
It remains submodular because adding the same vector commutes with
coordinatewise minimum and maximum. It also remains translation
invariant. Taking rectangular second differences and differentiating
translation invariance gives the Hessian identities
\[
 \partial_{ij}F_{d,\varepsilon}\le0\quad(i\ne j),
 \qquad
 \sum_j\partial_{ij}F_{d,\varepsilon}=0.
\]
Let $D=\Cov(Y)-\Cov(X)$ and
$\Delta_{ij}=D_{ii}+D_{jj}-2D_{ij}\ge0$.
For independent copies of $X,Y$, Gaussian interpolation and integration
by parts give
\[
 \frac d{dt}\E F_{d,\varepsilon}(\sqrt{1-t}\,X+\sqrt t\,Y)
 =-\frac14\sum_{i,j}\Delta_{ij}
       \E\,\partial_{ij}F_{d,\varepsilon}
             (\sqrt{1-t}\,X+\sqrt t\,Y)\ge0.
\]
The equality uses the zero row sums of the Hessian; the diagonal
terms vanish because $\Delta_{ii}=0$. More explicitly, the usual
interpolation derivative is
$\tfrac12\sum_{i,j}D_{ij}\E\partial_{ij}F_{d,\varepsilon}$,
which equals the displayed expression. Singular covariances cause no
problem: Gaussian integration by parts follows by writing each vector as
a linear image of a standard Gaussian vector, without requiring an
invertible matrix. Alternatively, add independent $N(0,\eta^2 I)$
noises to both vectors; this leaves all $\Delta_{ij}$ unchanged, and
Lipschitz continuity permits $\eta\downarrow0$ afterward.
The bounded derivatives justify interpolation, and the uniform
approximation permits $\varepsilon\downarrow0$.
\end{proof}

\begin{corollary}\label{cor:gaussian-paths}
Let $X,Y$ be centered Gaussian processes with continuous paths on $[0,A]$.
If
\[
 \E(X_t-X_s)^2\le\E(Y_t-Y_s)^2\qquad(s,t\in[0,A]),
\]
then, for every $d>0$,
\[
 \E\TV^d(X,[0,A])\le\E\TV^d(Y,[0,A]).
\]
\end{corollary}
\begin{proof}
Apply Lemma~\ref{lem:gaussian} on the dyadic grids
$\{jA/2^m:0\le j\le2^m\}$. Their subsequence truncated variations are
nondecreasing in $m$. Every fixed real-time partition can be approximated,
in order, by points of these grids. Continuity of each path therefore
makes the grid functionals increase to its full truncated variation.
Monotone convergence gives the expectation comparison, also when an
expectation is understood in the extended nonnegative reals.
\end{proof}

\section{The upper half: Brownian motion with Gaussian drift}\label{sec:upper}
\subsection{A tangent variogram and a Brownian motion with Gaussian drift}
Put $r=2H-1\in(0,1)$. Weighted arithmetic--geometric mean gives, for any
$T>0$ and $t\ge0$,
\[
 t^{1+r}\le (1-r)T^r t+rT^{r-1}t^2.
\]
Define
\[
 v=(1-r)T^r,\qquad w=rT^{r-1},\qquad
 Y_t=\sqrt v\,W_t+\sqrt w\,Zt,
\]
where $W$ is standard Brownian motion and $Z$ is an independent
standard normal. Thus $Y$ is a centered Gaussian process and
\[
 \E(B_t^H-B_s^H)^2\le v|t-s|+w(t-s)^2
 =\E(Y_t-Y_s)^2.
\]
The preceding lemma yields
$\E\TV^1(B^H,[0,L])\le\E\TV^1(Y,[0,L])$ for every $L$.

The following reflection calculation supplies a finite-horizon bound
with the Brownian drift rate of \cite{LochowskiMilos}; see also the
explicit truncated-variation formulas in \cite{LochowskiBM}.
Fix a drift $b\in\R$ and let $D_t=\sqrt v\,W_t+bt$. Consider its
reflection on $[-1/2,1/2]$:
\[
 X_t=X_0+D_t-K_t^++K_t^-,
\]
where $K^+,K^-$ are nondecreasing continuous processes starting at zero,
with their increments supported on $\{X=1/2\}$ and $\{X=-1/2\}$,
respectively. The invariant density is proportional to $e^{2bx/v}$:
it has zero flux $b\rho-(v/2)\rho'=0$. Start $X_0$ with this density,
independently of the future Brownian driver, so that $X$ is stationary.
Differentiating the logarithm of the density's normalizing integral gives
\[
 \E X=\tfrac12\coth(b/v)-\frac{v}{2b},
\]
with continuous value $0$ at $b=0$. It\^o's formula gives
\[
 X_A^2-X_0^2
 =2\sqrt v\int_0^A X_t\,dW_t+vA+2b\int_0^A X_t\,dt
   -K_A^+-K_A^-.
\]
The stochastic integral is square-integrable because $|X_t|\le1/2$;
the identity also shows integrability of $K_A^++K_A^-$. Stationarity
therefore implies
\[
 \E(K_A^++K_A^-)=A(v+2b\E X)=A b\coth(b/v),
\]
where the value at $b=0$ is $vA$.

Set $K=K^+-K^-$ and, crucially, $G=K-X_0$. Then
$D-G=X\in[-1/2,1/2]$, so $G$ is a feasible free-endpoint tube
approximation. The two regulator measures have disjoint supports, hence
$\TV(G,[0,A])=K_A^++K_A^-$. The tube identity yields the finite-horizon
estimate
\[
 \frac1A\E\TV^1(D,[0,A])\le b\coth(b/v).
\]
The auxiliary random initial value is harmless: for every realization,
$G$ is a feasible approximant to $D$. This estimate can be averaged
directly over $b=\sqrt w\,Z$, since
$0\le b\coth(b/v)\le v+|b|$. No interchange of expectation with a
long-time limit is required.
Combining it with the Gaussian comparison and
\eqref{eq:crossing-limit} gives
\begin{equation}\label{eq:gaussian-drift-upper}
 \boxed{c_H\le
 \E\!\left[\sqrt w\,Z\coth\!\left(\frac{\sqrt w\,Z}{v}\right)\right].}
\end{equation}

\subsection{A bound covering \texorpdfstring{$1/2<H\le3/4$}{1/2 < H <= 3/4}}
For $x\ge0$,
\[
 x\coth x\le1+x^2/3.
\]
For an elementary proof, differentiate
$(1+x^2/3)\sinh x-x\cosh x$:
its derivative is $x(x\cosh x-\sinh x)/3\ge0$.
The inequality is strict for $x>0$.
Equation \eqref{eq:gaussian-drift-upper} consequently gives
\[
 c_H\le v+\frac{w}{3v}
       =(1-r)T^r+\frac{r}{3(1-r)T}.
\]
The minimizing $T$ satisfies $T^{1+r}=1/[3(1-r)^2]$. Thus
\begin{equation}\label{eq:analytic-U}
 c_H\le U(H):=
 2H[2(1-H)]^{1/H-1}\,3^{-(2H-1)/(2H)}.
\end{equation}

We now check $U(H)<m_H$ for $1/2<H\le3/4$. Put $p=1/H\in[4/3,2)$
and let $R(p)=\log(U(1/p)/m_{1/p})$. Direct differentiation gives
\[
 R''(p)=\frac1{p(p-1)}
        -\frac14\psi_1\!\left(\frac{p+1}{2}\right)>0.
\]
Indeed, $p\le2$ implies the first term is at least $1/2$, whereas
$\psi_1((p+1)/2)\le\psi_1(1)=\pi^2/6<2$.
At $p=2$, $R(2)=0$. At $p=4/3$,
\[
 U(3/4)=\frac{3^{2/3}}{2^{4/3}},\qquad
 m_{3/4}=\frac{2^{2/3}\Gamma(7/6)}{\sqrt\pi},
\]
and $R(4/3)<0$ is equivalent to
\[
 \Gamma(7/6)>\frac{3^{2/3}\sqrt\pi}{4}.
\]
Here is a purely rational check of this last constant inequality.
The convergent expansion
\[
 \log\Gamma(1+x)=-\gamma x+
 \sum_{j\ge2}\frac{(-1)^j\zeta(j)}j x^j,\qquad |x|<1,
\]
follows by integrating \cite{DLMF}, equation 5.7.4. Its term magnitudes
decrease at $x=1/6$, so the alternating-series bound gives
\[
 \log\Gamma(7/6)>
 -\frac\gamma6+\frac{\pi^2}{432}-\frac{\zeta(3)}{648}
 >-\frac{1051}{12960},
\]
using $\gamma<3/5$, $\pi>3$, and
$\zeta(3)<1+1/8+\int_2^\infty t^{-3}\,dt=5/4$.
For completeness, strict midpoint convexity of $1/t$ gives
$\sum_{j=2}^n1/j<\log(n+1/2)-\log(3/2)$ and hence
$\gamma<1-\log(3/2)<3/5$. The last step uses
$\log(3/2)>2/5$, again by the strict midpoint rule. The same rule on
$[1,2]$ gives $\log2>2/3>\gamma$.
Since $1051/2160<487/1000$, the fifth alternating exponential
polynomial gives
\[
 \Gamma(7/6)^6>
 e^{-487/1000}>
 \sum_{k=0}^5\frac{(-487/1000)^k}{k!}
 =\frac{73733966284273793}{120000000000000000}
 >\frac{3072}{5000}.
\]
On the other hand, $\pi<22/7$ gives
\[
 \left(\frac{3^{2/3}\sqrt\pi}{4}\right)^6
 =\frac{81\pi^3}{4096}
 <\frac{81(22/7)^3}{4096}
 <\frac{307}{500}<\frac{3072}{5000}.
\]
This proves the endpoint inequality. Convexity of $R$ places its graph
below the chord between the two endpoints, so $R(p)<0$ for
$4/3\le p<2$. Hence \eqref{eq:analytic-U} proves the desired comparison
on $1/2<H\le3/4$.

\subsection{A bound covering \texorpdfstring{$3/4\le H<1$}{3/4 <= H < 1}}
Write $\mu=\E|Z|=\sqrt{2/\pi}$. For any $\sigma,v>0$, use
$\coth x=1+2/(e^{2x}-1)$ and the bound
$\varphi(z)\le\varphi(0)=\mu/2$ to obtain
\begin{align*}
 \E[\sigma Z\coth(\sigma Z/v)]
 &=\mu\sigma+
   4\sigma\int_0^\infty
      \frac{z\varphi(z)}{e^{2\sigma z/v}-1}\,dz\\
 &\le \mu\left(\sigma+\frac{\pi^2v^2}{12\sigma}\right).
\end{align*}
The integral evaluation uses
$\int_0^\infty u/(e^u-1)\,du=\sum_{n\ge1}n^{-2}=\pi^2/6$.

Choose $T=1$ in \eqref{eq:gaussian-drift-upper}. Then
$v=1-r$ and $\sigma=\sqrt r$. For $r\ge1/2$, put $x=\sqrt r$.
The inequality
\[
 x+\frac{\pi^2(1-x^2)^2}{12x}<1
\]
is equivalent, for $x<1$, to
\[
 \frac{\pi^2}{12}\frac{(1-x)(1+x)^2}{x}<1.
\]
The last rational function is
$x^{-1}+1-x-x^2$ and is decreasing. At $x=1/\sqrt2$ it equals
$1/2+1/\sqrt2$. The elementary estimates $\pi<22/7$ and
$\sqrt2>7/5$ give
\[
 \frac{\pi^2}{12}\left(\frac12+\frac1{\sqrt2}\right)
 <\frac{(22/7)^2}{12}\left(\frac12+\frac57\right)
 =\frac{8228}{8232}<1.
\]
Thus $c_H<\mu$ for every $3/4\le H<1$.

Finally, for $p>1$,
\[
 \frac d{dp}\log\E|Z|^p
 =\frac12\{\log2+\psi((p+1)/2)\}
 >\frac12(\log2-\gamma)>0.
\]
Therefore $m_H>\mu$ whenever $H<1$. This completes the proof of the
entire upper-half comparison.
\section{The lower half: periodic Brownian bridges}\label{sec:lower}

\subsection{A periodic Gaussian comparator}
Put $p=2H\in(0,1)$ and fix $T>0$. Define
\[
 v=(2-p)T^{p-1},\qquad w=(1-p)T^{p-2},\qquad
 L=\frac vw=\frac{2-p}{1-p}T.
\]
The function $t\mapsto t^{p-1}$ is convex. Its tangent at $T$ gives
\[
 t^p\ge vt-wt^2\qquad(t\ge0).
\]
For $t\ge0$, let $\tau(t)=t-L\lfloor t/L\rfloor$ and define
\[
 Y_t=\sqrt v\left(W_{\tau(t)}-\frac{\tau(t)}L W_L\right).
\]
Thus the same Brownian bridge of variance parameter $v$ and duration $L$
is repeated every period. The identified endpoint values are zero, so
$Y$ is a continuous centered Gaussian process. For phases
$\alpha,\beta\in[0,L]$, the bridge increment variance is
$v q(1-q/L)$, where $q=|\alpha-\beta|$. If $h=t-s\ge0$ and
$\delta=h\bmod L\in[0,L)$, then $q$ equals either $\delta$ or
$L-\delta$. Symmetry of this quadratic therefore gives
\[
 \mathbb E(Y_t-Y_s)^2=v\delta-w\delta^2
 \le\delta^p\le h^p.
\]
Thus fBm dominates this periodic process in increment variance on every
finite time interval, including increments that straddle period
boundaries or span several periods. Corollary~\ref{cor:gaussian-paths}
applies despite the repeated coordinates and resulting singular
covariance matrices.

\subsection{Cyclic crossing counts of a Brownian bridge}
Let $f$ be continuous on $[0,L]$ with $f(0)=f(L)$. For $c>0$, define
$N_c^\circ(f;a)$ as the largest integer $k\ge1$ for which there exist
\[
 0\le t_1\le\cdots\le t_{2k}\le L
\]
with alternating values below $a$ and above $a+c$: either
\[
 f(t_{2j-1})\le a,\qquad f(t_{2j})\ge a+c
 \quad(1\le j\le k),
\]
or the same conditions with odd and even indices interchanged. Set the
count to zero if there is no such sequence. Positive $c$ forces adjacent
witness times to be distinct. The first and last witnesses cannot be
both endpoints, because $f(0)=f(L)$ and they have opposite boundary
categories. Closing this witness word therefore certifies $k$ cyclic
upcrossings and $k$ cyclic downcrossings. Conversely, a circular
alternating boundary word gives such witnesses; thus this definition
agrees with the cyclic crossing count and is independent of the time cut.

Uniform continuity makes this count finite. It vanishes unless
$a\in[\min f,\max f-c]$. For each $k$, the set
$\{(f,a,c):c>0,\ N_c^\circ(f;a)\ge k\}$ is closed in the uniform path
topology and the parameter topology restricted to $c>0$: take convergent
witness sequences and extract a subsequence of their times from the
compact ordered simplex. Thus the count is jointly Borel measurable.
Its witness definition also shows directly that it is nonincreasing in
$c$ for fixed $f,a$. Define the measurable, finite path functional
\[
 \mathcal C_c(f)=2\int_{\mathbb R}N_c^\circ(f;a)\,da.
\]

\begin{lemma}[Expected cyclic crossing functional]\label{lem:cyclic}
For a Brownian bridge $Y$ of variance parameter $v$ and duration $L$,
\[
 \mathbb E\mathcal C_c(Y)
 =4c\sum_{n\ge1}n\exp\left(-\frac{2n^2c^2}{vL}\right).
\]
\end{lemma}

\begin{proof}
First take variance parameter one and duration $t$. Write
$p_t(x)=(2\pi t)^{-1/2}e^{-x^2/(2t)}$ and $R=\sqrt{2\lambda}$.
By translating a bridge pinned at zero, integration over strip locations
is the same as integration over starting points $x$ of bridges from $x$
to $x$, with the strip fixed at $[0,c]$. If $\PP^t_{x,x}$ denotes this
bridge law, put
\[
 A_n(t,c)=\int_\R\PP^t_{x,x}(N_c^\circ\ge n)\,dx.
\]
We calculate the Laplace transform in $t$ of $p_t(0)A_n(t,c)$; the
factor $p_t(0)$ is the endpoint density before conditioning.
The elementary Brownian formulas used below are
\[
 \E e^{-\lambda\tau_b}=e^{-bR}\quad(b>0),\qquad
 \int_0^\infty e^{-\lambda t}p_t(z)\,dt
 =R^{-1}e^{-R|z|},
\]
where $\tau_b$ is first passage through distance $b$. They follow from
the one-dimensional equation $u''/2=\lambda u$ and the heat resolvent.

If $x<0$, the event requires the usual $n$-th upcrossing, and successive
Brownian first-passage distances total $2nc-2x$, including the final heat
resolvent from the last upper-boundary hit back to $x$. Its transform is
\[
 R^{-1}e^{-(2nc-2x)R}.
\]
The case $x>c$ is symmetric. Roots at $0,c$ have zero measure in the
$x$ integral. Integrating the two exterior regions
gives $R^{-2}e^{-2ncR}$.

If $0<x<c$, let $K$ be the number of successive alternating boundary hits
in the rooted interval, counting whichever boundary is hit first and
ignoring repeated visits to the same side. The cyclic count is
$\lfloor K/2\rfloor$: when $K$ is even the time cut supplies the final
traversal, and when $K$ is odd the first and last sides agree. Therefore
$N_c^\circ\ge n$ means that the $2n$-th alternating hit occurs by time
$t$. The killed-interval exit transforms are
\[
 \mathbb E_x[e^{-\lambda\tau};B_\tau=0]
 =\frac{\sinh((c-x)R)}{\sinh(cR)},\qquad
 \mathbb E_x[e^{-\lambda\tau};B_\tau=c]
 =\frac{\sinh(xR)}{\sinh(cR)}.
\]
These formulas solve $u''/2=\lambda u$ with boundary values $(1,0)$ and
$(0,1)$, respectively. After the first hit, the remaining $2n-1$ strip traversals contribute
$e^{-(2n-1)cR}$; the final heat resolvent contributes
$R^{-1}e^{-(c-x)R}$ or $R^{-1}e^{-xR}$, respectively. The final
resolvent is unrestricted, since further traversals are permitted in
the event $N_c^\circ\ge n$. The interior integral is
\begin{align*}
 &\frac{e^{-(2n-1)cR}}{R\sinh(cR)}
 \int_0^c\left\{
 \sinh((c-x)R)e^{-(c-x)R}+\sinh(xR)e^{-xR}\right\}\,dx\\
 &\hspace{12mm}=
 e^{-2ncR}\left\{\frac{2c}{R(1-e^{-2cR})}-\frac1{R^2}\right\}.
\end{align*}
The exterior contribution cancels the last term. Hence
\[
 \int_0^\infty e^{-\lambda t}p_t(0)A_n(t,c)\,dt
 =\frac{2c}{R}\frac{e^{-2ncR}}{1-e^{-2cR}}.
\]
Summing the tail events over $n\ge1$ gives the transform of
$p_t(0)\int_\R\E N_c^\circ(Y;a)\,da$:
\[
 \frac{2c}{R}\frac{e^{-2cR}}{(1-e^{-2cR})^2}
 =\frac{2c}{R}\sum_{n\ge1}n e^{-2ncR}.
\]
All sums and integrals here are nonnegative, so Tonelli applies without
assuming integrability in advance. The computed transform is finite for
every $\lambda>0$, which gives local integrability in $t>0$. Since the
Laplace transform of $p_t(2nc)$ is $R^{-1}e^{-2ncR}$, uniqueness of
Laplace transforms and division by $p_t(0)$ give
\[
 \int_{\mathbb R}\mathbb E N_c^\circ(Y;a)\,da
 =2c\sum_{n\ge1}n e^{-2n^2c^2/t}.
\]
Multiplication by two gives the claimed crossing functional for almost
every positive duration. This extends to every duration without assuming
continuity of a hitting-time functional. By the witness definition,
$\mathcal C_c(f)$ is nonincreasing in $c$. Brownian scaling gives
$Q(t,c)/\sqrt t=Q(1,c/\sqrt t)$, where
$Q(t,c)=\mathbb E\mathcal C_c(Y)$ for a standard bridge of duration $t$.
Consequently $Q(t,c)/\sqrt t$ is nondecreasing in $t$, initially allowing
infinite values. The explicit series is finite and locally uniformly
convergent for $t>0$, hence continuous there. Squeezing between nearby
durations where equality holds proves finiteness and equality at every
duration.
Finally, variance parameter $v$ and duration $L$ are obtained by the time
change $t=vL$.
\end{proof}

\begin{lemma}[A deterministic lower bound]\label{lem:periodic-lower}
For a continuous $L$-periodic path $f$, an integer $n\ge1$, and $c>d>0$,
\[
 \TV^d(f,[0,nL])\ge n\mathcal C_c(f)-2\operatorname{range}(f).
\]
\end{lemma}
\begin{proof}
Take any approximant $u$ with $\|u-f\|_\infty\le d/2$. Choose a finite
partition of $[0,nL]$ on whose cells the oscillation of $f$ is less than
$(c-d)/2$, and let $g$ linearly interpolate the sampled values of $u$.
Then
\[
 \TV(g)\le\TV(u),\qquad
 \|g-f\|_\infty\le d/2+\max_{\text{cells}}\operatorname{osc}(f)<c/2.
\]
At a lower-boundary witness for the strip $[a,a+c]$, $g$ is strictly
below $y=a+c/2$; at an upper-boundary witness it is strictly above $y$.
For $N=N_c^\circ(f;a)\ge1$, concatenate a maximizing witness word over
the $n$ periods. This yields $2nN$ ordered alternating witnesses, hence
at least $2nN-1$, and in particular $2nN-2$, distinct crossings of $y$
by $g$ between consecutive witnesses.

For a continuous polygon $g$, let $J_g(y)$ count its level crossings.
Except at its finitely many vertex heights, this is the number of linear
segments with endpoint values strictly straddling $y$. Consequently the
elementary identity $|b-a|=\int\mathbf1_{\{\min(a,b)<y<\max(a,b)\}}\,dy$
gives
\[
 \TV(g)=\int_\R J_g(y)\,dy.
\]
Only strip positions in $[\min f,\max f-c]$ contribute, an interval
of length $(\operatorname{range}(f)-c)_+$. Integrating the forced crossing
count and using $dy=da$ therefore gives
\[
 \TV(u)\ge\TV(g)
 \ge n\mathcal C_c(f)-2(\operatorname{range}(f)-c)_+
 \ge n\mathcal C_c(f)-2\operatorname{range}(f).
\]
Taking the infimum in \eqref{eq:tube} proves the claim. This proof applies
even to discontinuous approximants $u$ and uses $c>d$ exactly to obtain
the strict interpolation margin.
\end{proof}

Apply Corollary~\ref{cor:gaussian-paths} and
Lemma~\ref{lem:periodic-lower} with $d=1$ and $c>1$:
\[
 \E\TV^1(B^H,[0,nL])
 \ge n\E\mathcal C_c(Y)-2\E\operatorname{range}(Y).
\]
The bridge range is integrable, since the bridge representation bounds
it by $2\sqrt v(\sup_{s\le L}|W_s|+|W_L|)$. Divide by $nL$ and use
\eqref{eq:crossing-limit} to obtain
\[
 c_H\ge\frac{4c}{L}\sum_{j\ge1}j
       \exp\left(-\frac{2j^2c^2}{vL}\right),\qquad c>1.
\]
The series converges locally uniformly in $c>0$, so letting $c\downarrow1$
gives the central lower estimate
\begin{equation}\label{eq:low-periodic}
 \boxed{c_H\ge\frac4L\sum_{n\ge1}n
       \exp\left(-\frac{2n^2}{vL}\right).}
\end{equation}

\subsection{A uniform elementary theta-sum bound}
\begin{lemma}
For every $a>0$,
\begin{equation}\label{eq:low-theta}
 \sum_{n\ge1}n e^{-an^2}\ge\frac1{2a}-\frac3{25}.
\end{equation}
\end{lemma}
\begin{proof}
For $0<a\le6/5$, apply Euler--Maclaurin through the $B_4$ term to
$f(x)=xe^{-ax^2}$ on $[0,\infty)$. It gives
\[
 \sum_{n\ge1}f(n)=\frac1{2a}-\frac1{12}-\frac a{120}+R,
 \qquad |R|\le\frac1{720}\int_0^\infty|f^{(4)}(x)|\,dx.
\]
The remainder bound uses $\sup_{0\le t\le1}|B_4(t)|=1/30$.
Writing $y=ax^2$,
\[
 f'''(x)=a(-6+24y-8y^2)e^{-y}.
\]
Its only positive-$x$ critical points correspond to
$y_\pm=(5\pm\sqrt{10})/2$. At either point its value divided by $a$ is
$(24-16y_\pm)e^{-y_\pm}$. The interior local maximum, after division
by $a$, is less than $4$: indeed
$y_->9/10$, $24-16y_-<48/5$, and
$e^{9/10}>12/5$, the latter following already from its Taylor polynomial
through degree three. The interior local minimum, after division by
$a$, is greater than $-1$: since $y_+>4$
and $(16y-24)e^{-y}$ decreases for $y\ge4$, its magnitude is less than
$40e^{-4}<1$. The inequality $e^4>40$ follows, for example, from the
Taylor polynomial through degree five. If $M$ and $m$ are these local
extremal values of $f'''$, its variation is
$6a+2M-2m<16a$, since it starts at $-6a$ and ends at zero. Thus
\[
 \frac1{2a}-\sum_{n\ge1}n e^{-an^2}
 <\frac1{12}+\frac a{120}+\frac a{45}
 =\frac1{12}+\frac{11a}{360}\le\frac3{25}.
\]

For $a\ge25/6$ the right side of \eqref{eq:low-theta} is nonpositive.
It remains to treat $6/5\le a\le25/6$, where the single term $e^{-a}$
suffices. Set
\[
 F(a)=e^a\left(\frac1{2a}-\frac3{25}\right).
\]
The sign of $F'$ is the sign of
$-6a^2+25a-25=-6(a-5/3)(a-5/2)$. Its maximum on the interval is therefore
attained at $6/5$ or $5/2$. The elementary Taylor bounds
\[
 e^{6/5}<\frac{583}{175}<\frac{10}{3},\qquad
 e^{5/2}<\frac{799}{64}<\frac{25}{2}
\]
give
\[
 F(6/5)<\frac{89}{90}<1,\qquad F(5/2)<1.
\]
For these rational Taylor upper bounds, truncate $e^x$ after degree $N$
and bound the remaining series by
$\frac{x^{N+1}}{(N+1)!}(1-x/(N+2))^{-1}$. Taking $N=2$ for $x=6/5$
and $N=3$ for $x=5/2$ gives the displayed fractions.
This proves \eqref{eq:low-theta} everywhere.
\end{proof}

\subsection{The interval \texorpdfstring{$1/3\le H<1/2$}{1/3 <= H < 1/2}}
Put $r=1-2H\in(0,1/3]$ and $k=3/25$. Since
$v=(1+r)T^{-r}$ and $L=(1+r)T/r$, the theta bound gives
\[
 c_H\ge v-\frac{4k}{L}
 =(1+r)T^{-r}-\frac{4kr}{(1+r)T}.
\]
Choose
\[
 T=\left(\frac{4k}{(1+r)^2}\right)^{1/(1-r)}.
\]
Then
\begin{equation}\label{eq:low-central}
 c_H\ge\mathcal L(r):=
 (1-r)(4k)^{-r/(1-r)}(1+r)^{(1+r)/(1-r)}.
\end{equation}
We show that $\mathcal L(r)>m_H$ on this whole interval. Since
\[
 m_H=\frac{2^{1/(1-r)}\Gamma(1/2+1/(1-r))}{\sqrt\pi},
\]
differentiation gives
\[
 \frac d{dr}\log\frac{\mathcal L(r)}{m_H}
 =\frac{N(r)}{(1-r)^2},\qquad
 N(r)=-\log(8k)+2\log(1+r)
       -\psi\left(\frac12+\frac1{1-r}\right).
\]
At zero,
\[
 N(0)=\gamma-2+\log(25/6)>0.
\]
Here is an entirely rational check of the positive margin. The convex
trapezoid rule for $1/x$, summed over $[n,\infty)$, gives
$H_n-\log n-\gamma<1/(2n)$, where $H_n$ is the $n$-th harmonic number.
At $n=5$ this yields $\gamma>131/60-\log5$. The alternating Taylor
bound for $\log(1+1/5)$ now gives
\[
 N(0)>\frac{11}{60}-\log(6/5)
 >\frac{11}{60}-\left(\frac15-\frac1{50}+\frac1{375}\right)
 =\frac1{1500}>0.
\]

The convex midpoint rule gives, for $z>1/2$,
\[
 \psi_1(z)=\sum_{j\ge0}\frac1{(z+j)^2}
 \le\int_{-1/2}^\infty\frac{dx}{(z+x)^2}
 =\frac1{z-1/2}.
\]
Consequently
\[
 N'(r)=\frac2{1+r}
 -\frac{\psi_1(1/2+1/(1-r))}{(1-r)^2}
 \ge\frac2{1+r}-\frac1{1-r}
 =\frac{1-3r}{1-r^2}\ge0.
\]
The logarithmic ratio is zero at $r=0$ and has strictly positive
derivative for $0<r\le1/3$. Hence $c_H\ge\mathcal L(r)>m_H$ for
$1/3\le H<1/2$.

\subsection{The interval \texorpdfstring{$0<H\le1/3$}{0 < H <= 1/3}}
Set $x=1/(2H)\ge3/2$. Retain only the $n=1$ term in
\eqref{eq:low-periodic} and choose the scale so that $2/(vL)=x$.
Equivalently,
\[
 T=\left[\frac{2(x-1)}{(2x-1)^2}\right]^x.
\]
It follows that
\[
 c_H\ge\frac{4(x-1)}{2x-1}
 \left[\frac{(2x-1)^2}{2e(x-1)}\right]^x.
\]
The centered Stirling (Burnside) bound
\[
 \Gamma(x+1/2)<\sqrt{2\pi}(x/e)^x\qquad(x>0)
\]
therefore yields
\[
 \frac{c_H}{m_H}>
 R(x):=\frac{2\sqrt2(x-1)}{2x-1}
 \left[1+\frac1{4x(x-1)}\right]^x.
\]
The rational terms cancel on differentiation:
\[
 \frac d{dx}\log R(x)
 =\log\left[1+\frac1{4x(x-1)}\right]>0.
\]
Since $R(3/2)=8/(3\sqrt6)>1$, this proves $c_H>m_H$ for every
$0<H\le1/3$ as well.

For completeness, the Gamma inequality follows by setting
$g(x)=\log\Gamma(x+1/2)-x\log x+x-\tfrac12\log(2\pi)$.
The digamma Laplace representation \cite{DLMF}, equation 5.9.12,
combined with the elementary Frullani integral for $\log x$, gives
\[
 g'(x)=\int_0^\infty e^{-xt}
 \left\{\frac1t-\frac1{2\sinh(t/2)}\right\}\,dt>0,
\]
The integral converges for every $x>0$: its bracket is $O(t)$ at zero
and $O(1/t)$ at infinity. Stirling's formula \cite{DLMF}, equation
5.11.1, gives $g(x)\to0$ at infinity. Hence $g(x)<0$.

\section{Completion and scope}
Section~\ref{sec:upper} proves the strict upper comparison for
$1/2<H\le3/4$ and $3/4\le H<1$. Section~\ref{sec:lower} proves the
strict lower comparison for $0<H\le1/3$ and $1/3\le H<1/2$.

The Brownian equality also follows from the estimates already proved.
Reflection with $v=1$ and $b=0$ gives $c_{1/2}\le1$. Conversely,
repeat a standard Brownian bridge of any duration $L>0$. Its increment
variance is $\delta-\delta^2/L\le|t-s|$, so the same Gaussian and
cyclic comparisons, followed by \eqref{eq:low-theta}, give
\[
 c_{1/2}\ge\frac4L\sum_{n\ge1}n e^{-2n^2/L}
 \ge1-\frac{12}{25L}.
\]
Letting $L\to\infty$ shows $c_{1/2}=1=m_{1/2}$. This completes the
proof of Theorem~\ref{thm:main}.

The proof concerns the comparison with the Gaussian moment; it does not
provide an exact formula for $c_H$ at non-Brownian parameters. In
particular, neither the Gaussian-drift upper bound nor the periodic-bridge
lower bound is asserted to equal $c_H$.

\section*{Acknowledgments}
The authors acknowledge the assistance of OpenAI's ChatGPT in developing
proof strategies, deriving and checking intermediate estimates, verifying
references, and preparing the exposition. Responsibility for the
mathematical content and the final manuscript rests with the authors.

\end{document}